\documentclass[final]{siamltex}
\newtheorem{thm}{Theorem}[section]

\newtheorem{Lemma}{Observasion}[section]
\newtheorem{Def}{Definition}[section]
\title{Considering Adelman's Shortest Permutation strings}
\author{Hesam Dashti\thanks{Department of Mathematics, University of Wisconsin, Madison, USA}}
\begin{document}
\maketitle
\begin{abstract}
In this report, we consider Adelman's algorithm for generating shortest permutation strings. We introduce a new representation approach which reveals some properties of Adelman's algorithm. 
\end{abstract}
%\begin{AMS}
%15A15, 15A09, 15A23
%\end{AMS}
\pagestyle{myheadings}
\thispagestyle{plain}
\markboth{Hesam Dashti}{Considering Adelman's Shortest Permutation strings}
\section{Introduction}
Since Donald Knuth in 1971 \cite{Knuth} published a set of open problems with computational flavor, the problem of generating a permutation string with the least length has attracted a great deal of attention. For the past four decades, some algorithms have been introduced for tackling this problem but still many ambiguities remain. In addition to the beauty of the problem, its applicability in the security area motivated researchers to focus on this combinatorics problem.\\
Newey \cite{Newey} and Adelman \cite{Adel} proved an upper bound (n$^2$-2n+4) for the length of a $\Pi_n$-Complete word. Later on some other proofs have been proposed and confirmed the boundary. Newey determined the lower bounds for n$\leq$7 as shown in Table \ref{T1}. In order to have a better understanding of the definition of $\Pi_n$-Complete strings here we note the following examples:\\
121 is a $\Pi_2$-Complete word.\\
1221 is a $\Pi_2$-Complete word.\\
1231213 is a $\Pi_3$-Complete word.\vspace{2ex}\\
\begin{table}[h]\label{T1}
\caption{The length of permutation strings for n$\leq$7, proved by Newey}
\begin{center}
\begin{tabular}{|l|l|l|l|l|l|l|l|}
\hline n & 1 & 2 & 3 & 4 & 5 & 6 & 7\\
\hline Lower bound & 1 & 3 & 7 & 12 & 19 & 28 & 39\\
\hline
\end{tabular}
\end{center}
\end{table}\\
After introducing the open problem of the lower bound of $\Pi_n$-Complete string, the next section focused on needed preliminaries for revisiting the Adelman's Algorithm. 
\section{Preliminaries}\label{Prelim}
Since our algorithm is based on Adelman's Algorithm \cite{Adel}, the reader is strongly encouraged to read Adelman's paper first.\\
Let us start with introducing our notation:\\
\textbf{Notation:}
\begin{itemize}
\item $S_n$=\{1, $\ldots$, n\} is a set of alphabets.
\item $W\in S^*_n$ is a sentence that may, or may not, be a $\Pi_n$-Complete string.
\item A Permutation Word (PW) is a permutation of $S$ without repetition.
\item $\delta$ is a generic symbol of PW's.
\item x-permutations is a generic symbol for all the permutations of length x, where x $\leq$ n.
\end{itemize}
\subsection{Revisiting Adelman's Algorithm}\label{RevAdel}
In this section we recall Adelman's Algorithm \cite{Adel} and introduce a framework to capture its permutation strings. We will use this framework for the rest of the paper.\vspace{2ex}\\
\textbf{Adelman's Algorithm:}\\
\begin{Def}\label{TPrime}
A string $W=w_1w_2\ldots,w_k$ is a R$_n$-String iff satisfies the following conditions:\\
(a) $W \in S^*_n$\\
(b) $w_{i+1}$=$w_i$(mod n)+1
\end{Def}\vspace{2ex}\\
Example: The followings are R$_7$-String:\\
1234\\
234567123456712345\vspace{2ex}\\
Here, we use this definition to construct a $\Pi_n$-Complete string.\\
\begin{itemize}
\item[i] Construct a unique string T$_n \in S^*_{n-1}$ such that
\begin{itemize}
\item[a)] T$_n$ is a R$_{n-1}$-String,
\item[b)] T$_n$ is of length n$^2$-3n+4,
\item[c)] The first letter of T$_n$ is 1.
\end{itemize}
\item[ii] Construct T'$_n$ = T$_n$ and then change it as follows: for all i, if 1$\leq$i$\leq$n-2, then insert the letter n into T'$_n$ after the i'th occurrence of the letter n$-i$.
\item[iii] Construct $\Pi_n$ = nT'$_n$n.
\end{itemize}
Here we borrow two examples from the Adelman's paper which are depicted in Table \ref{AdelAlgEx}.
\begin{table}[h!]
\caption{The three steps of constructing $\Pi_4$ and $\Pi_7$ are shown}
\begin{center}
\begin{tabular}{|l|l|l|}
\hline & $\Pi_4$ & $\Pi_7$\\
\hline T & 12312312 & 12345612345612345612345612345612\\
\hline T' & 1234124312 & 1234567123457612347561237456127345612\\
\hline $\Pi_n$ & 412341243124 & 712345671234576123475612374561273456127\\
\hline
\end{tabular}
\end{center}
\label{AdelAlgEx}
\end{table}\\
Thereafter, Adelman proved the $\Pi_n$ is a $\Pi_n$-Complete string of length n$^2$-2n+4 which is described in details in his paper \cite{Adel}. In order to rewrite these complete strings, we use a sliding window of length $\|S_n\|$, which traces the string ($W$) and extracts substrings of length $n$. Each substring would be a PW of $S$ if we duplicate the last letter of a substring at the beginning of the next substring. For example, assume we want to rewrite the above $\Pi_4$-Complete in this framework: 412341243124 = 4123 4124 3124. Now, if we perform the duplication process we would have 412341243124 = 4123 \underline{3}412 \underline{2}431 \underline{1}24. The latter looks perfect since we have exactly 3 PW without any extra letters. Let us do the same on $\Pi_7$ and compare the output; $\Pi_7$ = 7123456 712345\underline{7}6 1234756 1237456 1273456 127. As it is shown, in addition to the extra letters we have a repeated letter in the second PW which changes the attitude and this substring will not be a PW. Now, let us check the string with duplications: $\Pi_7$ = 7123456 \underline{6}712345 \underline{5}761234 \underline{4}756123 \underline{3}745612 \underline{2}734561 \underline{1}27. This is a set of perfect PW's with two extra letters at the end. Using this representation we would be able to generate all the different $\Pi_n$-Complete strings when the letters have been substituted. For example, we can consider another $\Pi_7$-Complete string where 7 is substituted by 2: $\Pi_7$ = 2173456 \underline{6}217345 \underline{5}261734 \underline{4}256173 \underline{3}245617 \underline{7}234561 \underline{1}72 which is a $\Pi_7$-Complete string.\\
To conclude this section, we can say that the new representation of the Adelman's Algorithm gives a structure for $\Pi_n$-Complete strings. For example, we derived the template of $\Pi_7$-Complete strings as shown in the Table ~\ref{ExamAdel1}. 
\begin{table}[h!]
\caption{The structure and some examples of the $\Pi_7$ based on Adelman's Algorithm. As it is shown, the examples exactly follow the structure with assigning different letters to the structure elements ($a_{i}'s$).}
\begin{center}
\begin{tabular}{|l|l|l|}
\hline Structure & Ex1 & Ex2 \\
\hline $a_1a_2a_3a_4a_5a_6a_7 $ & 1234567 & 7123456\\
\hline $a_7a_1a_2a_3a_4a_5a_6 $ & 7123456 & 6712345 \\
\hline $a_6a_1a_7a_2a_3a_4a_5$ & 6172345 & 5761234 \\
\hline $a_5a_1a_6a_7a_2a_3a_4 $ & 5167234 & 4756123\\
\hline $a_4a_1a_5a_6a_7a_2a_3 $ & 4156723 & 3745612\\
\hline $a_3a_1a_4a_5a_6a_7a_2 $ & 3145672 & 2734561 \\
\hline $a_2a_3a_1$ & 231 & 127\\
\hline
\end{tabular}
\end{center}
\label{ExamAdel1}
\end{table}\\
We experimentally verified the completeness of Adelman's Algorithm as described in the next section. Based on these investigations we have extracted the following observations:\\
\begin{Lemma}\label{Lem1}
For a given set $S_n$=\{1, $\ldots$, n\}, Adelman's rule generates a string $W$ which covers n-permutations. Note that this rule does not mean that a substring $w \in W$ of length k$^2$-2k+4 (k$<$n) can cover every k-permutations of words S$_k \in$ S$_n$.
\end{Lemma}\\
\emph{Example:}\\
912345678\\
891234567\\
798123456\\
679\\
This string does not cover every 4-Permutations. 
\begin{thm}\label{aolast}
For a string generated by Adelman's rule, substitution of the last two letters does not effect the completeness of the string. 
\end{thm}\\
\emph{Example:}\\
{\small W1= A123456789\\
W2= 9A12345678\\
W3= 8A91234567\\
W4= 7A89123456\\
W5= 6A78912345\\
W6= 5A67891234\\
W7= 4A56789123\\
W8= 3A45678912\\
W9= 2A34567891\\
W10= 12A (or 1A2)}\vspace{2ex}\\
\textbf{Proof}\\
In order to proof this theorem we consider different positions of 'A' and '2' with respect to each other.\\
When $\delta_9$ =2, the '2' in the W9 would be selected and we should have a complete PW after it. Hence, it does not matter what the arrangement of 'A' and '2' is in W10. Similarly, when the $\delta_{10}$ = 2, the last '2' would be selected where there is a complete PW between the '2' in W10 and W9 without considering the last 'A' in W10. \\
Hence, position of 'A' in W10 does not affect the position of '2'.\\
When $\delta_9$ = 'A', as shown before the 'A' from W9 would be chosen. So the position of '2' in W10 is not important. In a case that $\delta_{10}$ = 'A', we would choose the last 'A' to fill the $\delta$, so the string before that would be:\\ 
{\small W1= 123456789\\
W2= 912345678\\
W3= 891234567\\
W4= 789123456\\
W5= 678912345\\
W6= 567891234\\
W7= 456789123\\
W8= 345678912\\
W9= 234567891\\
W10= 1 (or 12)\\}
In both cases, we have more than 9 PWs for the 9 letters of $\delta$. So, arrangement of the last line is not important.\\
\begin{Lemma}\label{lem2}
For a given set, $S_n$=\{$a_1, \ldots, a_n$\}, and a string $W \in S$ (of length n$^2$-2n+4) generated according to Adelman's rule, a string $W'$ ($W' \in S, W' \subset W$) could be modified to a $\Pi_k$-Complete by removing $n-k$ letters of $W'$, such that $W'$ satisfies Adelman's rule. 
\end{Lemma}\vspace{2ex}\\
Actually, this Observation is based on the fact that removing 'n-k' letters (excluding the first letter of $W'$) from $W'$ yields to a $\Pi_k$-Complete string of $S'_k$=\{$b_1, \ldots, b_k$\} $\in$ $S_n$=\{$a_n, \ldots, a_n$\}.\vspace{1ex}\\
\emph{Example:}\\
In the previous example we saw the following string is not a $\Pi_4$-Complete string.\\
912345678\\
891234567\\
798123456\\
679\\
Here, by removing every 5 letters we would have a $\Pi_4$-Complete string. Example:\\
9128\\
8912\\
9812\\
9\\
which is the same as:\\
9128\\
8912\\
2981\\
129\\
This string is a $\Pi_4$-Complete string. Note that by removing every 5 letters (except 9) the remaining string follows from Adelman's algorithm and is a complete string. 
\begin{Lemma}\label{Lem3}
There exists a slight modification of Adelman's Algorithm that generates a $\Pi_n$-Complete string of length n$^2$-2n+5. This modification changes the cyclic behavior of the letters, in one permutation word by repeating the first letter at the end of the permutation word.
\end{Lemma}\vspace{2ex}\\
\emph{Example:} An example is provided in the (Table \ref{T7}-W). In this example, W6 and W7 end with a$_4$ and rearranging this string gives us W' as it shown in (Table \ref{T7}-W').
\begin{table}[h!]
\caption{A $\Pi_9$-Complete string with a small modification on the circular behavior in Adelman's Algorithm}
\begin{center}
\begin{tabular}{|l|l|l|}
\hline & W & W'\\
\hline W1 & $ a_1a_2a_3a_4a_5a_6a_7a_8a_9$ & $a_1a_2a_3a_4a_5a_6a_7a_8a_9$\\
\hline W2 & $ a_9a_1a_2a_3a_4a_5a_6a_7a_8$& $a_9a_1a_2a_3a_4a_5a_6a_7a_8$\\
\hline W3 & $a_8a_1a_9a_2a_3a_4a_5a_6a_7$ & $a_8a_1a_9a_2a_3a_4a_5a_6a_7$\\
\hline W4 & $a_7a_1a_8a_9a_2a_3a_4a_5a_6$& $a_7a_1a_8a_9a_2a_3a_4a_5a_6$\\
\hline W5 & $a_6a_1a_7a_8a_9a_2a_3a_4a_5$ & $a_6a_1a_7a_8a_9a_2a_3a_4a_5$\\ 
\hline W6 & $a_5a_1a_6a_7a_8a_9a_2a_3a_4$& $a_5a_1a_6a_7a_8a_9a_2a_3a_4$\\
\hline W7 & $a_1a_5a_6a_7a_8a_9a_2a_3a_4$& $a_4a_1a_5a_6a_7a_8a_9a_2a_3a_4$\\
\hline W8 & $a_4a_1a_5a_6a_7a_8a_9a_2a_3$& $a_4a_1a_5a_6a_7a_8a_9a_2a_3$ \\
\hline W9 & $a_3a_4a_1$& $a_3a_4a_1$\\
\hline
\end{tabular}
\end{center}
\label{T7}
\end{table}\\
This string follows Adelman's rule, except for W7 that ends with its starting letter (a$_4$) and also, W8 does not have cyclic behavior. This string is a $\Pi_9$-Complete string and the reader can verify this as shown in the next section.
\section{Diving into Adelman's Algorithm}\label{AppA}
Let us start with an example of Adelman's Algorithm for an alphabet set $S_{10}$=\{1, $\ldots$, 9, A\}. Based on Adelman's rule, the $W$ (Table \ref{T4}) is a $\Pi_{10}$-Complete string. It is clear that the ending letter of each line (a PW) is repeated at the beginning of the next line and should be ignored in terms of counting the length. These substrings are shown in Table\ref{T4}. 
Since the theoretical terminology of Adelman's Algorithm is proved in his paper, here we empirically see why this string covers all the permutations $\delta$ of $S$. To do so, we consider different $\delta$'s with different positions of a letter $\alpha \in S$.\\
\textbf{Considering permutation words subject to positions of the letter 'A' in $\delta$:}\\
In each of the following parts of this section, we fix position of the letter 'A' in the string $\delta$ and find the best 'A' from the string 'W'; such that substrings before and after this 'A' in 'W' can cover every necessary permutation of substrings before and after the 'A' in $\delta$.\\
\underline{$\delta_1$='A'}: Since 'A' is the first letter of the strings $\delta$ and $W$, the rest of $W$ includes 9 complete permutation words, (plus \underline{two extra letters}) which can be used to generate every permutation words of length 9. Hence, all the permutations, $\delta$, starting with 'A' would be covered.\\
\underline{$\delta_2$='A'}: Since there is a PW before the second 'A', the second 'A' can be used for filling the $\delta_2$ position (Table \ref{T4}). Each line is a PW, so, we need to fill 8 letters of $\delta$ and we have 8 PW which means all of the 8-permutations would be covered.\\
\underline{$\delta_3$='A'}: As shown in the Table \ref{T4} we filled 3 letters and need 7 more letters, where we have 7 PW's. \\
The rest of positions are similar, just let us consider the last case:\\
\underline{$\delta_{10}$='A'}: This means we need to fill 9 letters, from the PW's before the last 'A' in 'W'. We have 9 PWs for the 9 blank positions, so we would find all of the permutations.\\ 
\begin{table}[h!]
\caption{Considering different positions of the letter 'A' in a $\delta$}
\begin{center}
\begin{tabular}{|l|l|l|l|l|l|}
\hline & W & $\delta_1$ = A & $\delta_2$=A & $\delta_3$=A & $\delta_{10}=A$\\
\hline W1 & A123456789 & 123456789 && &123456789\\
\hline W2 & 9A12345678 & 912345678 &12345678 & &912345678\\
\hline W3 & 8A91234567 & 891234567 &891234567 & 891234567 & 891234567\\
\hline W4 & 7A89123456 & 789123456 &789123456 & 789123456&789123456\\
\hline W5 & 6A78912345 & 678912345 &678912345 & 678912345&678912345\\ 
\hline W6 & 5A67891234 & 567891234 &567891234 & 567891234&567891234\\
\hline W7 & 4A56789123 & 456789123 &456789123 & 456789123&456789123\\
\hline W8 & 3A45678912 & 345678912 &345678912 & 345678912&345678912\\
\hline W9 & 2A34567891 & 234567891 &234567891 & 234567891&234567891\\
\hline W10 & 12A & 12 & 12 & 12 & 12\\
\hline
\end{tabular}
\end{center}
\label{T4}
\end{table}\vspace{2ex}\\
Now, we perform a similar process for the letter '9'. Strings of different steps are shown in the associated columns in Table\ref{T5}.\\
\underline{$\delta_{1}$='9'}: Using the first '9' in $W$, we want to find 9 letters from the remaining string (Table \ref{T5}). In a recursive manner from the Adelman's rule, we can generate all the 8-permutations from this string.\\
\underline{$\delta_{2}$='9'}: Since we have a complete PW before the first '9', we have the previous sequence but this time to cover 8 letters. With this string we can cover every 9-permutations, but need 8-permutation words, \underline{so we have one extra PW}.\\
\underline{$\delta_{3}$='9'}: Since we covered the first two positions of '9' by the first '9', in this case, we use the second '9'. \\
Before the second '9' we have two PW's (red colored), so we can cover every 2-permutations. For the remaining positions of $\delta$ we have blue colored strings (Table \ref{T5}). By ignoring the last 7 letters in W3, the rest covers every 7-permutations. Hence, we have a \underline{7 extra letters} at the beginning. Other cases are similar to the previous ones, let us jump to the last case:\\
\underline{$\delta_{10}$='9'} Here, we are interested in finding 9 letters from the string before the last '9' in $W$ as depicted in the last column of Table \ref{T5}. The first line is a PW and covers at least one letter. The rest covers all 8-permutations, where the last 5 letters are unnecessary. Hence, we cover every 10-permutations ending with '9', where there are \underline{5 extra letters}.
\begin{table}[h]
\caption{Considering different positions of the letter '9' in a $\delta$}
\begin{center}
\begin{tabular}{|l|l|l|l|l|l|}
\hline & W & $\delta_1$=9 & $\delta_{2}$='9' & $\delta_{3}$='9'&$\delta_{10}$='9'\\
\hline W1 & A123456789 & & & \underline{A123456789}&A12345678\\
\hline W2 & 9A12345678 & A12345678 & A12345678& \underline{9A12345678}& A12345678\\
\hline W3 & 8A91234567 & 8A1234567 & 8A1234567& \underline{8A}1234567&8A1234567\\
\hline W4 & 7A89123456 & 7A8123456 & 7A8123456& 7A8123456&7A8123456 \\
\hline W5 & 6A78912345 & 6A7812345 & 6A7812345& 6A7812345&6A7812345\\ 
\hline W6 & 5A67891234 & 5A6781234 & 5A6781234& 5A6781234&5A6781234\\
\hline W7 & 4A56789123 & 4A5678123 & 4A5678123& 4A5678123&4A5678123\\
\hline W8 & 3A45678912 & 3A4567812 & 3A4567812& 3A4567812&3A4567812\\
\hline W9 & 2A34567891 & 2A3456781 & 2A3456781& 2A3456781&2A345678\\
\hline W10 & 12A & 12A & 12A & 12 A & \\ 
\hline
\end{tabular}
\end{center}
\label{T5}
\end{table}\\
In the following we consider another letter ('4') and the rest of letters and positions are similar.\\
\underline{$\delta_{1}$='4'}: the first '4' would be chosen, the remaining string is shown in Table \ref{T6}. Ignoring the first line, the rest covers all 9-permutations. Hence, in this case we have \underline{5 extra letters}.\\
\underline{$\delta_{2}$='4'}: And again \underline{4 extra letters}. Let us jump to the last case:\\
\underline{$\delta_{10}$='4'}: Based on Adelman's rule, this sequence covers all 9-permutations.
\begin{table}[h]\label{T6}
\caption{Considering different positions of the letter '4' in a $\delta$}
\begin{center}
\begin{tabular}{|l|l|l|l|l|l|}
\hline & W & $\delta_1$='4' & $\delta_{2}$='4' & $\delta_{10}$='4'\\
\hline W1 & A123456789 & 56789& & A12356789\\
\hline W2 & 9A12345678 & 9A1235678& 5678& 9A1235678\\
\hline W3 & 8A91234567 & 8A9123567& 8A9123567& 8A9123567\\
\hline W4 & 7A89123456 & 7A8912356& 7A8912356& 7A8912356\\
\hline W5 & 6A78912345 & 6A7891235& 6A7891235& 6A7891235\\ 
\hline W6 & 5A67891234 & 5A6789123 & 5A6789123& 5A6789123\\
\hline W7 & 4A56789123 & A56789123 & A56789123& A56789123\\
\hline W8 & 3A45678912 & 3A5678912& 3A5678912& 3A5678912 \\
\hline W9 & 2A34567891 & 2A3567891& 2A3567891& 2A3\\
\hline W10 & 12A & 12A & 12A & \\ 
\hline
\end{tabular}
\end{center}
\end{table}\\

\end{document}